\catcode`\@=11

\font\teneufm=eufm10
\font\seveneufm=eufm7
\font\fiveeufm=eufm5
\newfam\eufmfam
\textfont\eufmfam=\teneufm  \scriptfont\eufmfam=\seveneufm
     \scriptscriptfont\eufmfam=\fiveeufm
\def\eufm@{\hexnumber@\eufmfam}

\font\tenmsa=msam10
\font\sevenmsa=msam7
\font\fivemsa=msam5
\font\tenmsb=msbm10
\font\sevenmsb=msbm7
\font\fivemsb=msbm5
\newfam\msafam
\newfam\msbfam
\textfont\msafam=\tenmsa  \scriptfont\msafam=\sevenmsa
     \scriptscriptfont\msafam=\fivemsa
\textfont\msbfam=\tenmsb  \scriptfont\msbfam=\sevenmsb
     \scriptscriptfont\msbfam=\fivemsb

\def\hexnumber@#1{\ifnum#1<10 \number#1\else
    \ifnum#1=10 A\else\ifnum#1=11 B\else\ifnum#1=12 C\else
    \ifnum#1=13 D\else\ifnum#1=14 E\else\ifnum#1=15
F\fi\fi\fi\fi\fi\fi\fi}

\def\msa@{\hexnumber@\msafam}
\def\msb@{\hexnumber@\msbfam}

\mathchardef\gx="2\eufm@78
\mathchardef\gf="2\eufm@66
\mathchardef\gg="2\eufm@67
\mathchardef\gh="2\eufm@68
\mathchardef\gm="2\eufm@6D
\mathchardef\gp="2\eufm@70
\mathchardef\gd="2\eufm@64
\mathchardef\ga="2\eufm@61

\mathchardef\boxdot="2\msa@00
\mathchardef\boxplus="2\msa@01
\mathchardef\boxtimes="2\msa@02
\mathchardef\square="0\msa@03
\mathchardef\blacksquare="0\msa@04
\mathchardef\centerdot="2\msa@05
\mathchardef\lozenge="0\msa@06
\mathchardef\blacklozenge="0\msa@07
\mathchardef\circlearrowright="3\msa@08
\mathchardef\circlearrowleft="3\msa@09
\mathchardef\rightleftharpoons="3\msa@0A
\mathchardef\leftrightharpoons="3\msa@0B
\mathchardef\boxminus="2\msa@0C
\mathchardef\Vdash="3\msa@0D
\mathchardef\Vvdash="3\msa@0E
\mathchardef\vDash="3\msa@0F
\mathchardef\twoheadrightarrow="3\msa@10
\mathchardef\twoheadleftarrow="3\msa@11
\mathchardef\leftleftarrows="3\msa@12
\mathchardef\rightrightarrows="3\msa@13
\mathchardef\upuparrows="3\msa@14
\mathchardef\downdownarrows="3\msa@15
\mathchardef\upharpoonright="3\msa@16

\mathchardef\downharpoonright="3\msa@17
\mathchardef\upharpoonleft="3\msa@18
\mathchardef\downharpoonleft="3\msa@19
\mathchardef\rightarrowtail="3\msa@1A
\mathchardef\leftarrowtail="3\msa@1B
\mathchardef\leftrightarrows="3\msa@1C
\mathchardef\rightleftarrows="3\msa@1D
\mathchardef\Lsh="3\msa@1E
\mathchardef\Rsh="3\msa@1F
\mathchardef\rightsquigarrow="3\msa@20
\mathchardef\leftrightsquigarrow="3\msa@21
\mathchardef\looparrowleft="3\msa@22
\mathchardef\looparrowright="3\msa@23
\mathchardef\circeq="3\msa@24
\mathchardef\succsim="3\msa@25
\mathchardef\gtrsim="3\msa@26
\mathchardef\gtrapprox="3\msa@27
\mathchardef\multimap="3\msa@28
\mathchardef\therefore="3\msa@29
\mathchardef\because="3\msa@2A
\mathchardef\doteqdot="3\msa@2B

\mathchardef\triangleq="3\msa@2C
\mathchardef\precsim="3\msa@2D
\mathchardef\lesssim="3\msa@2E
\mathchardef\lessapprox="3\msa@2F
\mathchardef\eqslantless="3\msa@30
\mathchardef\eqslantgtr="3\msa@31
\mathchardef\curlyeqprec="3\msa@32
\mathchardef\curlyeqsucc="3\msa@33
\mathchardef\preccurlyeq="3\msa@34
\mathchardef\leqq="3\msa@35
\mathchardef\leqslant="3\msa@36
\mathchardef\lessgtr="3\msa@37
\mathchardef\backprime="0\msa@38
\mathchardef\risingdotseq="3\msa@3A
\mathchardef\fallingdotseq="3\msa@3B
\mathchardef\succcurlyeq="3\msa@3C
\mathchardef\geqq="3\msa@3D
\mathchardef\geqslant="3\msa@3E
\mathchardef\gtrless="3\msa@3F
\mathchardef\sqsubset="3\msa@40
\mathchardef\sqsupset="3\msa@41
\mathchardef\trianglerighteq="3\msa@44
\mathchardef\trianglelefteq="3\msa@45
\mathchardef\bigstar="0\msa@46
\mathchardef\between="3\msa@47
\mathchardef\blacktriangledown="0\msa@48
\mathchardef\blacktriangleright="3\msa@49
\mathchardef\blacktriangleleft="3\msa@4A
\mathchardef\blacktriangle="0\msa@4E
\mathchardef\triangledown="0\msa@4F
\mathchardef\eqcirc="3\msa@50
\mathchardef\lesseqgtr="3\msa@51
\mathchardef\gtreqless="3\msa@52
\mathchardef\lesseqqgtr="3\msa@53
\mathchardef\gtreqqless="3\msa@54
\mathchardef\Rrightarrow="3\msa@56
\mathchardef\Lleftarrow="3\msa@57
\mathchardef\veebar="2\msa@59
\mathchardef\barwedge="2\msa@5A
\mathchardef\doublebarwedge="2\msa@5B
\mathchardef\angle="0\msa@5C
\mathchardef\measuredangle="0\msa@5D
\mathchardef\sphericalangle="0\msa@5E
\mathchardef\varpropto="3\msa@5F
\mathchardef\smallsmile="3\msa@60
\mathchardef\smallfrown="3\msa@61
\mathchardef\Subset="3\msa@62
\mathchardef\Supset="3\msa@63
\mathchardef\Cup="2\msa@64

\mathchardef\Cap="2\msa@65

\mathchardef\curlywedge="2\msa@66
\mathchardef\curlyvee="2\msa@67
\mathchardef\leftthreetimes="2\msa@68
\mathchardef\rightthreetimes="2\msa@69
\mathchardef\subseteqq="3\msa@6A
\mathchardef\supseteqq="3\msa@6B
\mathchardef\bumpeq="3\msa@6C
\mathchardef\Bumpeq="3\msa@6D
\mathchardef\lll="3\msa@6E

\mathchardef\ggg="3\msa@6F

\mathchardef\circledS="0\msa@73
\mathchardef\pitchfork="3\msa@74
\mathchardef\dotplus="2\msa@75
\mathchardef\backsim="3\msa@76
\mathchardef\backsimeq="3\msa@77
\mathchardef\complement="0\msa@7B
\mathchardef\intercal="2\msa@7C
\mathchardef\circledcirc="2\msa@7D
\mathchardef\circledast="2\msa@7E
\mathchardef\circleddash="2\msa@7F
\def\ulcorner{\delimiter"4\msa@70\msa@70 }
\def\urcorner{\delimiter"5\msa@71\msa@71 }
\def\llcorner{\delimiter"4\msa@78\msa@78 }
\def\lrcorner{\delimiter"5\msa@79\msa@79 }
\def\yen{\mathhexbox\msa@55 }
\def\checkmark{\mathhexbox\msa@58 }
\def\circledR{\mathhexbox\msa@72 }
\def\maltese{\mathhexbox\msa@7A }
\mathchardef\lvertneqq="3\msb@00
\mathchardef\gvertneqq="3\msb@01
\mathchardef\nleq="3\msb@02
\mathchardef\ngeq="3\msb@03
\mathchardef\nless="3\msb@04
\mathchardef\ngtr="3\msb@05
\mathchardef\nprec="3\msb@06
\mathchardef\nsucc="3\msb@07
\mathchardef\lneqq="3\msb@08
\mathchardef\gneqq="3\msb@09
\mathchardef\nleqslant="3\msb@0A
\mathchardef\ngeqslant="3\msb@0B
\mathchardef\lneq="3\msb@0C
\mathchardef\gneq="3\msb@0D
\mathchardef\npreceq="3\msb@0E
\mathchardef\nsucceq="3\msb@0F
\mathchardef\precnsim="3\msb@10
\mathchardef\succnsim="3\msb@11
\mathchardef\lnsim="3\msb@12
\mathchardef\gnsim="3\msb@13
\mathchardef\nleqq="3\msb@14
\mathchardef\ngeqq="3\msb@15
\mathchardef\precneqq="3\msb@16
\mathchardef\succneqq="3\msb@17
\mathchardef\precnapprox="3\msb@18
\mathchardef\succnapprox="3\msb@19
\mathchardef\lnapprox="3\msb@1A
\mathchardef\gnapprox="3\msb@1B
\mathchardef\nsim="3\msb@1C
\mathchardef\napprox="3\msb@1D
\mathchardef\nsubseteqq="3\msb@22
\mathchardef\nsupseteqq="3\msb@23
\mathchardef\subsetneqq="3\msb@24
\mathchardef\supsetneqq="3\msb@25
\mathchardef\subsetneq="3\msb@28
\mathchardef\supsetneq="3\msb@29
\mathchardef\nsubseteq="3\msb@2A
\mathchardef\nsupseteq="3\msb@2B
\mathchardef\nparallel="3\msb@2C
\mathchardef\nmid="3\msb@2D
\mathchardef\nshortmid="3\msb@2E
\mathchardef\nshortparallel="3\msb@2F
\mathchardef\nvdash="3\msb@30
\mathchardef\nVdash="3\msb@31
\mathchardef\nvDash="3\msb@32
\mathchardef\nVDash="3\msb@33
\mathchardef\ntrianglerighteq="3\msb@34
\mathchardef\ntrianglelefteq="3\msb@35
\mathchardef\ntriangleleft="3\msb@36
\mathchardef\ntriangleright="3\msb@37
\mathchardef\nleftarrow="3\msb@38
\mathchardef\nrightarrow="3\msb@39
\mathchardef\nLeftarrow="3\msb@3A
\mathchardef\nRightarrow="3\msb@3B
\mathchardef\nLeftrightarrow="3\msb@3C
\mathchardef\nleftrightarrow="3\msb@3D
\mathchardef\divideontimes="2\msb@3E
\mathchardef\varnothing="0\msb@3F
\mathchardef\nexists="0\msb@40
\mathchardef\mho="0\msb@66
\mathchardef\thorn="0\msb@67
\mathchardef\beth="0\msb@69
\mathchardef\gimel="0\msb@6A
\mathchardef\daleth="0\msb@6B
\mathchardef\lessdot="3\msb@6C
\mathchardef\gtrdot="3\msb@6D
\mathchardef\ltimes="2\msb@6E
\mathchardef\rtimes="2\msb@6F
\mathchardef\shortmid="3\msb@70
\mathchardef\shortparallel="3\msb@71
\mathchardef\smallsetminus="2\msb@72
\mathchardef\thicksim="3\msb@73
\mathchardef\thickapprox="3\msb@74
\mathchardef\approxeq="3\msb@75
\mathchardef\succapprox="3\msb@76
\mathchardef\precapprox="3\msb@77
\mathchardef\curvearrowleft="3\msb@78
\mathchardef\curvearrowright="3\msb@79
\mathchardef\digamma="0\msb@7A
\mathchardef\varkappa="0\msb@7B
\mathchardef\hslash="0\msb@7D
\mathchardef\hbar="0\msb@7E
\mathchardef\backepsilon="3\msb@7F
\def\Bbb{\ifmmode\let\next\Bbb@\else
    \def\next{\errmessage{Use \string\Bbb\space only in math
mode}}\fi\next}
\def\Bbb@#1{{\Bbb@@{#1}}}
\def\Bbb@@#1{\fam\msbfam#1}

\catcode`\@=\active

\def\eps{{\epsilon}}

\def\<{\langle}
\def\>{\rangle}

\def\tens{\mathop{\otimes}}
\def\ra{{\triangleleft}}

\def\Ad{{\rm Ad}}

\def\ev{{\rm ev}}

\def\id{{\rm id}}

\def\sw#1{{\sb{(#1)}}}
\def\su#1{{\sb{[#1]}}}
\def\ot{\otimes}

\def\extd{{\rm d}}
\def\cS{\mathcal{S}}

\def\proof{\goodbreak\noindent{\bf Proof\quad}}

\def\text#1{{\rm #1}}
\def\note#1{}

\documentclass[12pt]{article}
\newtheorem{lemma}{Lemma}[section]
\newtheorem{propos}[lemma]{Proposition}
\newtheorem{example}[lemma]{Example}
\newtheorem{theorem}[lemma]{Theorem}
\newtheorem{cor}[lemma]{Corollary}

\newtheorem{defin}[lemma]{Definition}
\newtheorem{remark}[lemma]{Remark}

\begin{document}



\hfill Swan-Maths-03/18

\begin{center} {\Large The van Est spectral sequence for Hopf
algebras}
\\ \baselineskip 13pt{\ }
{\ }\\  E. J. Beggs \& Tomasz Brzezi\'nski \\
{\ } \\ Department of Mathematics\\
University of Wales, Swansea\\ Wales SA2 8PP
\end{center}

\begin{quote}
\noindent{\bf Abstract.}
Various aspects of the de Rham cohomology of Hopf algebras are
discussed.
In particular, it is shown that the de Rham cohomology of an algebra
with the
differentiable coaction of a cosemisimple Hopf algebra with trivial
0-th cohomology
group, reduces to the de Rham cohomology of (co)invariant forms.
Spectral
sequences are discussed and the van Est spectral sequence for Hopf
algebras
is introduced. A definition of
  Hopf-Lie algebra cohomology is also given.
\end{quote}

\section{Introduction}
The idea of calculating the de-Rham cohomology of a manifold direcly
from
the definition is a rather daunting one.  Fortunately there are all
sorts of tools, from the Meyer Vietoris exact sequence to the
equivalence of many types of cohomology given by sheaf theory, to
assist.  The situation in noncommutative geometry is not so fortunate
in this regard.  In this paper we will give some results on the de
Rham cohomology of Hopf algebras equipped with bicovariant
differential calculi \cite{worondiff}.  These results will be the
noncommutative geometry analogues of well known results about Lie
groups.

The paper \cite{CheEil} considered the connection between the
de Rham cohomology of a Lie group and the cohomology of the left
invariant differential forms.  In the simplest case, the cohomologies
coincide for a compact connected Lie group, with compactness
appearing in the
guise of a normalised Haar integral.  The cohomology of the left
invariant differential forms can be taken to define the Lie algebra
cohomology \cite{JacLie}.  The van Est spectral sequence \cite{BorWal}
is a more general result, whose statement requires the additional
concept of group cohomology \cite{BroCoh}.  The involvement of the Lie
algebra cohomology is fortunate in that for a (usual finite
dimensional) Lie algebra finding it only requires a finite dimensional
computation.

In this paper we give results for Hopf algebras
(with bicovariant differential calculus) corresponding to the results
for
Lie groups. First we show that if a Hopf algebra has a normalised
left integral, then
its de Rham cohomology is isomorphic to the cohomology of left
invariant forms.
Then we prove a version of the van Est spectral sequence.
It then remains to identify the cohomology of left invariant forms
in an easily calculable form. This must be, by definition, the
cohomology of the Hopf-Lie algebra of the vector fields on the
Hopf algebra with given
bicovariant differential structure. For Hopf-Lie or braided Lie algebras
see \cite{worondiff,brLie}.
  The usual definition of
Lie algebra cohomology must be modified, as the Lie bracket no longer
satisfies the Jacobi identities.

Throughout the paper we work with vector spaces, algebras etc.\ over
a commutative field $k$. The unadorned tensor product is always over
$k$. Given an algebra $A$, $\Omega A = \bigoplus_{n\geq 0} \Omega^nA$
denotes a graded differential algebra with the differential $\extd :
\Omega^nA\to \Omega^{n+1}A$ such that $A = \Omega^0A$. Each of the
$\Omega^A$ is an $A$-bimodule and we always assume that  $\Omega A$
satisfies the {\em density condition}, i.e., that  for all $\omega\in
\Omega^1A$ there exist $a_i, b_i\in A$, $i=1,2,\ldots , m$ such that
$\omega = \sum_{i=1}^M a_i.\extd b_i$. We often refer to $\Omega A$
as a {\em differential calculus} or a {\em differential structure on
$A$}. For any  Hopf algebra $P$, the coproduct is denoted by
$\Delta$, the counit by $\eps$ and the antipode by $S$. We use the
Sweedler sigma-notation for a coproduct (without a sigma), i.e., we
write
$\Delta(p) = p\sw 1\ot p\sw 2$ (summation understood). A left
$P$-coaction
on $F$, $\lambda:F\to P\ot F$ is denoted on elements by the Sweedler
notation
as well, but with indices in square brackets, i.e., $\lambda(f)  =
f\su{-1}\ot f\su 0$
(summation understood). The coassociativity of $\lambda$ entails that
for all
$f\in F$, $f\su{-1}\sw 1 \ot f\su{-1}\sw 2\ot f\su 0= f\su{-1} \ot
f\su{0}\su{- 1}\ot
f\su 0\su 0$, thus we simply write  $f\su{-2}\ot f\su{-1}\ot f\su 0$
etc. For the
theory of bicovariant differential calculi on a Hopf algebra we refer
to \cite{worondiff} or to \cite{KliSch:gro}.

\section{Coactions on the de Rham cohomology}\label{lppp}
In this section we suppose that $P$ is a Hopf algebra with a
bicovariant differential calculus, and that $P$ left coacts on an
algebra $M$ (with a given differential structure)
    by a differentiable map $\lambda:M\to P\tens M$.
   The differentiability just means that $\lambda$
   extends to a map
$\lambda_*:\Omega^n M \to
\Omega^n(P\tens M)$ which commutes with the differential $\extd$
(i.e.\ a cochain map).

\begin{remark}
\rm  By
definition of the
tensor product differential structure,
\[
\Omega^{n}(P\tens M)\,=\,(\Omega^{0}P\tens \Omega^{n}M) \oplus
(\Omega^{1}P\tens \Omega^{n-1}M) \oplus\dots\oplus
(\Omega^{n}P\tens \Omega^{0}M)\ ,
\]
and the differential $\extd$ on $\Omega^{n}(P\tens M)$ corresponds to
$\extd\tens\id+(-1)^{r}\id\tens \extd$ on $\Omega^{r}P\tens
\Omega^{s}M$.
We define projections $\Pi_{r}:\Omega^{n}(P\tens M)\to
\Omega^{r}P\tens \Omega^{n-r}M$.  Then there is a left $P$-coaction on
$\Omega^n M$ given by $\bar\lambda=\Pi_0\circ\lambda_*:\Omega^n M\to
P\tens \Omega^n
M$.
By definition of $\extd$ on the tensor product,
$\Pi_0(\extd\xi)=(\id\tens\extd)\Pi_0(\xi)$
for $\xi\in\Omega^*(P\tens M)$ and from this we see that
$\extd:\Omega^n M\to \Omega^{n+1} M$
is a left $P$-comodule map. It follows that there is a left
$P$-coaction
$\tilde\lambda:H^*_{dR}(M)\to P\tens H^*_{dR}(M)$
    given by $[\omega]\mapsto (\id\tens[\bullet ])\bar\lambda(\omega)$.
\end{remark}

The next proposition could be regarded as a part of a noncommutative
K\"unneth theorem.

\begin{propos} \label{ccoo}   The image of the
left coaction $\tilde\lambda:H_{dR}^{n}(M)\to P\tens
H_{dR}^{n}(M)$ is contained in $(\ker \extd:P\to\Omega^{1}P)\tens
H_{dR}^{n}(M)$.
\end{propos}
\proof  Given $\omega\in \Omega^{n}M$ with
$\extd \omega=0$, we set
$\Pi_{0}(\omega)=\sum_{i}p_{i}\tens\tau_{i}\in
P\tens\Omega^{n}M$ and
$\Pi_{1}(\omega)=\sum_{j}\xi_{j}\tens\eta_{j}\in
\Omega^{1}P\tens\Omega^{n-1}M$.  Since $\extd\omega=0$ we have
$0=\sum_{i}p_{i}\tens\extd\tau_{i}\in P\tens\Omega^{n+1}M$ and
$0=\sum_{i}\extd p_{i}\tens\tau_{i}-\sum_{j}\xi_{j}\tens\extd
\eta_{j}\in
\Omega^{1}P\tens\Omega^{n}M$.  Without loss of generality, from the
first equality we may assume that all $\tau_{i}\in\ker \extd:
\Omega^{n}M\to \Omega^{n+1}M$.  Now using the quotient map
$[\bullet ]:(\ker \extd:
\Omega^{n}M\to \Omega^{n+1}M)\to H_{dR}^{n}(M)$, the second equation
gives
$\sum_{i}\extd p_{i}\tens[\tau_{i}]=0\in \Omega^{1}P\tens
H_{dR}^{n}(M)$.
It follows that
$\sum_{i} p_{i}\tens[\tau_{i}]\in (\ker \extd:P\to\Omega^{1}P)\tens
H_{dR}^{n}(M)$.\quad$\square$

\begin{defin}\label{connop}\rm
    Let $P$ be a Hopf algebra with a given bicovariant differential
structure. $P$  is
called a {\em connected Hopf algebra} if
$H_{dR}^0(P)=k$. The unit element $1\in k$  is identified with the
class of the identity $1_P$ in
$H_{dR}^0(P)$.
\end{defin}

Note that the notion of connectedness introduced here is
{\em differential calculus dependent},
i.e., a Hopf algebra can be  a connected Hopf algebra with respect to
a given
differential structure and does not have be a connected Hopf algebra
with
respect to another differential structure.
For example the quantum group $SU_q(2)$ is a connected Hopf algebra
with
respect to the 4D-differential calculi of Woronowicz (cf.\
\cite{worondiff}). The quantum group $GL_q(2)$ is not connected with
respect to these calculi, as the quantum determinant induces a
non-trivial class in the de Rham cohomology. On the other hand, any
Hopf algebra (over a field) is a connected Hopf algebra with respect
to the
universal differential structure.

\begin{cor}
    If $P$ is a connected Hopf algebra (with its given bicovariant
differential structure), then
all elements of $H^*_{dR}(M)$ are fixed by the coaction
$\tilde\lambda$.
\end{cor}

In \cite{RemBic} it was shown that $H^0_{dR}(P)$ is a Hopf algebra.
We have shown that
$([\bullet ]\tens\id)\tilde\lambda:H^*_{dR}(M)\to H^0_{dR}(P)\tens
H^*_{dR}(M)$ is
a left $H^0_{dR}(P)$-coaction. We can go further to a coaction of the
entire
graded Hopf algebra $H^*_{dR}(P)$, also described
in
\cite{RemBic}.  To do this we will state the K\"unneth theorem for
noncommutetive de Rham cohomology.  Its proof is standard and
straightforward homological algebra, but it is useful to state it in
this context.

\begin{theorem}
    Let $M$ and $N$ be algebras with differential calculi, and give
$N\tens M$ the standard
tensor product differential calculus. Then there is an isomorphism
\[
\bigoplus_{n\ge r\ge 0} H_{dR}^r(N) \hat{\tens} H_{dR}^{n-r}(M) \cong
H_{dR}^n(N\tens M)
\]
given by mapping $[\omega]\hat{\tens}[\xi]\in H_{dR}^r(N) \hat{\tens}
H_{dR}^{n-r}(M)$
    to $[\omega\tens \xi]\in H_{dR}^n(N\tens M)$. The $\hat{\tens}$
operation is
   the standard tensor product in which the wedge product
becomes modified by the grading
to give $(x\hat{\tens} y)\wedge (w\hat{\tens} z)=(-1)^{nm}\, (x\wedge
w) \hat{\tens} (y\wedge z)$,
where $y\in H_{dR}^{n}(M)$ and $w\in H_{dR}^m(N)$.

In addition, if there are differentiable algebra maps $\phi:N\to N'$
and
$\psi:M\to M'$, then in terms of the isomorphism above we have
$\phi_*\hat{\tens}\psi_*:
H_{dR}^r(N) \hat{\tens} H_{dR}^{s}(M)\to H_{dR}^r(N') \hat{\tens}
H_{dR}^{s}(M')$
corresponding to $(\phi\tens\psi)_*:H_{dR}^n(N\tens M)\to
H_{dR}^n(N'\tens M')$.
\end{theorem}

\begin{cor}
     The graded Hopf algebra $H^*_{dR}(P)$ coacts on
$H^*_{dR}(M)$ by
\[
H^n_{dR}(P) \stackrel{\lambda_*}{\longrightarrow} H^n_{dR}(P\tens
M)\,\cong\,
\bigoplus_{n\ge r\ge 0} H_{dR}^r(P) \hat{\tens} H_{dR}^{n-r}(M)\ .
\]
\end{cor}

\section{Integrals and invariant forms}
Again we suppose that $P$ is a Hopf algebra with a
bicovariant differential calculus, and that $P$ left coacts on an
algebra $M$ (with a given differential structure)
    by a differentiable map $\lambda:M\to P\tens M$.  Recall the
    definition of a normalised left integral on a Hopf algebra.

\begin{defin} \rm A {\em left integral} on a Hopf algebra $P$ is a
linear map $\int:P\to k$ such that
       $(\int\tens\id)\Delta=I_P.\int:P\to P$. A left integral is said
to be {\em normalised} provided
   $\int I_p=1$.
\end{defin}

Throughout this section we suppose that $P$ has a normalised left
integral $\int$. Since we are working over a field this is equivalent
to assuming that
$P$ is a {\em cosemisimple Hopf algebra} (i.e., a sum of simple
coalgebras)   \cite[14.0.3]{Swe:Hop}. Given a left
$P$-comodule $E$ with coaction $\bar\lambda$, we define a map
$$
\Bbb I=(\int\tens\id)\bar\lambda:E \to {}^{coP} E:=\{\omega\in E \;
|\; \bar{\lambda}(\omega) = 1_P\otimes \omega\}.
$$

\begin{lemma} The left invariant forms ${}^{coP}(\Omega^n M) := \{
\omega\in \Omega^n M \; |\; \bar{\lambda}(\omega) = 1_P\otimes
\omega\}$ form a cochain
complex with the
usual de Rham differential. Furthermore, the map $\Bbb I:\Omega^n M
\to
{}^{coP}(\Omega^n M)$ is a
cochain map.
\end{lemma}
\proof As $\extd:\Omega^n M\to \Omega^{n+1} M$
is a left $P$-comodule map it follows that $\extd$ preserves the
invariant forms.  Also we find that
\[
\Bbb I(\extd\omega)\,=\,(\int\tens\id)\bar\lambda(\extd\omega)
\,=\, (\int\tens\id)(\id\tens\extd)\bar\lambda(\omega)\,=\,\extd (\Bbb
I(\omega))\ .\quad\square
\]

\medskip
Now we have two cochain maps, $\Bbb I:\Omega^n M \to
{}^{coP}(\Omega^n M)$
and the inclusion map $i: {}^{coP}(\Omega^n M) \to\Omega^n M $. As
the
integral is normalised it follows
that $ \Bbb I\circ i:{}^{coP}(\Omega^n M) \to {}^{coP}(\Omega^n M)$
is the
identity. Next consider the
induced cohomology maps $H(\Bbb I):H^n_{dR}(M)\to
H^n({}^{coP}(\Omega^*
M),\extd)$ and
$H(i): H^n({}^{coP}(\Omega^* M),\extd)\to H^n_{dR}(M)$. We see that
$H(\Bbb I)\circ H(i)$ is the identity on $H^n({}^{coP}(\Omega^*
M),\extd)$, and as a result
$H(i)\circ H(\Bbb I)$ is a projection on $H^n_{dR}(M)$. It is obvious
that
the image of $H(i)$ is contained in ${}^{coP}(H^n_{dR}(M))$, but it
is
less obvious that
the image is ${}^{coP}(H^n_{dR}(M))$.

\begin{propos} The image of $H(i): H^n({}^{coP}(\Omega^* M),\extd)\to
H^n_{dR}(M)$
is ${}^{coP}(H^n_{dR}(M))$.
\end{propos}
\proof Suppose that $[\omega]\in H^n_{dR}(M)$ and that
$(\id\tens[\bullet ])\bar\lambda(\omega)=1_P\tens[\omega]$.
Write $\bar\lambda(\omega)=\sum_i p_i\tens\tau_i$, where $\extd
\tau_i=0$. Then
$\sum_i p_i\tens[\tau_i]=1_P\tens [\omega]$, or
\[
\sum_i p_i\tens\tau_i\,=\,1_P\tens\omega\,+\,\sum_j q_j\tens
\extd\eta_j\ .
\]
Applying $\int\tens\id$ to this we get
\[
\Bbb I(\omega)\,=\,\sum_i \int(p_i)\,\tau_i\,=\,\omega\,+\,\sum_j
\int(q_j)\, \extd\eta_j\ ,
\]
showing that $[\Bbb I(\omega)]=[\omega]\in H^n_{dR}(M)$, and
$[\Bbb I(\omega)]$ is in the image of $H(i)$. \quad$\square$
\medskip

The above discussion leads to the main result of this section.

\begin{theorem}\label{thm.iioo} Let $P$ be a Hopf algebra with a
bicovariant
differential calculus. Suppose that
$P$ coacts on $M$ by a differentiable left coaction.

\medskip \noindent
1)\quad If  $P$ has a normalised left integral, then there is an
isomorphism $H(i): H^n({}^{coP}(\Omega^*
M),\extd)\to
{}^{coP}(H^n_{dR}(M))$.

\smallskip \noindent
2)\quad If $P$ is connected then ${}^{coP}(H^n_{dR}(M))=H^n_{dR}(M)$.

\smallskip \noindent
3)\quad If $P$ is a connected cosemisimple Hopf algebra, then
${}^{coP}(H^n_{dR}(M))=H^n_{dR}(M)= H^n({}^{coP}(\Omega^*
M),\extd)$.

\end{theorem}

\begin{cor} \label{iioo} If $P$ is a Hopf algebra with bicovariant
differential structure,
and $P$ has a normalised left integral, then the inclusion map
induces an isomorphism $H(i): H^n({}^{coP}(\Omega^* P),\extd)\to
{}^{coP}(H^n_{dR}(P))$.
If $P$ is connected then ${}^{coP}(H^n_{dR}(P))=H^n_{dR}(P)$.
\end{cor}

\begin{remark}\rm
In \cite[Theorem~3.1]{drhodquant} the cohomology of
Hopf algebras $SL_q(N)$, $SO_q(N)$ and $Sp_q(N)$ was
shown by explicit calculation
to be the same as the cohomology of their left invariant forms for
certain bicovariant differential calculi. All
these Hopf algebras are cosemisimple, i.e., they
have a normalised left integral, and they can be shown to be
connected,
hence this part of \cite[Theorem~3.1]{drhodquant}
follows from
Corollary~\ref{iioo}. On the other hand it is also computed in
\cite[Theorem~3.1]{drhodquant}
that for $O_q(2n+1)$  and certain bicovariant differential calculi,
the de Rham
cohomology decomposes into a sum of two copies of the de Rham
cohomology of
invariant forms. Although $O_q(2n+1)$  are cosemisimple Hopf
algebras, they are
not connected with respect to differential calculi discussed in
\cite{drhodquant}. In
addition to the de Rham cohomology class $1_P$ there is another
non-trivial class
induced by the quantum determinant. Although
\cite[Theorem~3.1]{drhodquant}
cannot be inferred from Theorem~\ref{thm.iioo} and
Corollary~\ref{iioo}, the latter
indicate the origin of this decomposition: every `connected component'
leads to
contributions to the de Rham cohomology that go beyond the invariant
part.
Similarly, the decomposition of the de
Rham cohomology of $GL_q(N)$ in \cite[Theorem~3.1(2)]{drhodquant} can
be
expected from Theorem~\ref{thm.iioo} in view of the fact that the
powers of the
quantum determinant induce nontrivial cohomology classes in
$H^0_{dR}(GL_q(N))$.
\end{remark}

\section{Hopf cochain cohomology}
In this section $P$ is a Hopf algebra.

\begin{defin} \label{yuyu} \rm Suppose that $F$ is a left
$P$-comodule.  Define $D^{n}=P^{\tens n+1}\tens F$ for $n\ge 0$, with
the tensor product left coaction.  The map $\extd:D^{n}\to D^{n+1}$ is
defined by
       \[
       \extd(p_{0}\tens\ldots\tens p_{n}\tens f)\,=\, \sum_{n+1\ge i\ge
       0}(-1)^{i}\, p_{0}\tens\ldots\tens p_{i-1}\tens 1 \tens
       p_{i}\tens\ldots \tens p_{n}\tens f\ .
       \]
       It follows that $\extd$ is left $P$-covariant, and that
$d^{2}=0$.
       The cohomology of the $P$-invariant
       complex $({}^{coP}D^{n},\extd)$ is called a {\em Hopf cochain
       cohomology of $P$ with coefficients in $F$} and is denoted by
       $H_{c}^{n}(P;F)$.
\end{defin}

The reader familiar with the cohomology theory of algebras (and  with
the descent theory
in particular) will recognise in $(D^n,\extd)$ in 
Definition~\ref{yuyu} the {\em Amitsur complex} of the algebra $P$
\cite{Ami:sim} (cf.\ \cite[Section~6]{Art:Azu}). Note, however, that, motivated by 
the group cohomology, for a Hopf cochain cohomology we take the 
coinvariant part of the Amitsur complex.

\begin{propos}\label{yooy}
     Let $F$ be a left $P$-comodule with coaction $\lambda: F\to P\ot
F$. Suppose that there is a left action $\mu :P\tens F\to F$
      which is a left $P$-comodule map with respect to $\lambda$ and
the tensor
      product coaction, i.e., such that, for all $p\in P$ and $f\in
F$,
$\lambda(\mu(p\otimes f)) = p\sw 1f\su{-1}\otimes \mu(p\sw 2\ot f\su
0)$.
Then $H_{c}^{n}(P;F)=0$ for $n\ge 1$
and $H_{c}^{0}(P;F)\cong {}^{coP} F$,
with the isomorphism mapping $f\in {}^{coP} F$ to $[1_P\tens f]\in
H_{c}^{0}(P;F)$.
\end{propos}
\proof Define the map $h:P^{\tens n+2}\tens F =D^{n+1} \to P^{\tens
n+1}\tens
F=D^n$ for $n\ge -1$ by
\[
h(p_{0}\tens\ldots \tens p_{n+1}\tens f)\,=\,
(-1)^{n+1}\,p_{0}\tens\ldots \tens
p_{n}\tens p_{n+1}.  f\ .
\]
Now we calculate
\begin{eqnarray*}
       h\,\extd(p_{0}\tens\ldots\tens p_{n+1}\tens f) &=& \sum_{n+1\ge
i\ge
       0}(-1)^{i+n}\, p_{0}\tens\ldots\tens p_{i-1}\tens 1 \tens
       p_{i}\tens\ldots \tens p_{n+1}. f  \cr
&&+\, p_{0}\tens\ldots \tens p_{n+1}\tens f\ ,\cr
\extd h(p_{0}\tens\ldots \tens p_{n+1}\tens f) &=& (-1)^{n+1}\,1
\tens p_{0}\tens\ldots \tens
p_{n}\tens p_{n+1}.  f + \dots \cr
&& + \,
p_{0}\tens\ldots \tens
p_{n}\tens 1\tens p_{n+1}.  f\ .
\end{eqnarray*}
   From this we see that
\[
(d\circ h+h\circ d)(p_{0}\tens\ldots \tens p_{n+1}\tens f)\,=\,
p_{0}\tens\ldots \tens p_{n+1}\tens f\ .
\]
Since $h$ is left $P$-covariant, we have a cochain homotopy which
contracts the complex,
showing that $H_{c}^{n}(P;F)=0$ for $n\ge 1$. To find $H_{c}^{0}(P;F)$
we need to be
more careful, as our hypothetical $D^{-1}$ is not part of the complex.
However
we can still write $x=\extd h(x)+h(\extd x)$ for $x\in D^0$. Then if
$x\in {}^{coP} D^0$
and $dx=0$, we have $x=1\tens h(x)$, so we can identify $\ker
\extd:D^0\to D^1$
as $1_P\tens {}^{coP} F$. \quad $\square$
\medskip

A left $F$-comodule $P$  that satisfies assumptions of
Proposition~\ref{yooy}, i.e., such that it is also a left $P$-module
with an
action that is compatible with a coaction is known as  a left
{\em Hopf module} (cf.\ \cite{Swe:Hop}). As an explicit example one
can
consider {\em cleft extensions} of algebras
(cf.\ \cite{BlaCoh:cro}, \cite{DoiTak:cle}).
\begin{example}\label{ex.cleft}\rm
Let $F$ be a left $P$-comodule algebra with the coaction $\lambda$,
and
let $M= {}^{coP}F$ be the algebra of coinvariants. Suppose there
exists a
left $P$-colinear map $\Phi: P\to F$ such that $\Phi(1) = 1$ and
$\Phi$ is
convolution-invertible, i.e., there is a map $\Phi^{-1}:P\to F$ such
that for
all $p\in P$, $\Phi(p\sw1)\Phi^{-1}(p\sw 2) =
\Phi^{-1}(p\sw1)\Phi(p\sw 2) = \eps(p)1_P$.
Then $F$ is called a {\em cleft extension of $M$}. This is an example
of a
Hopf-Galois extension that has a geometric meaning of a trivial
principal
bundle (cf.\ \cite{BrzMaj:gau}). One can show that $F\cong P\otimes
M$ as
a right $M$-module and as a left $P$-comodule, where the coaction in
$P\otimes M$ is given by $\Delta\otimes \id$ (cf.\
\cite[Theorem~9]{DoiTak:cle}).
Explicitly, the isomorphism $\Theta: F\to P\otimes M$ and its inverse
are
$$
\Theta(f) = f\su{-2}\ot \Phi^{-1}(f\su{-1})f\su 0, \qquad
\Theta^{-1}(p\otimes x) = \Phi(p)x.
$$
The  map  $\Theta$ is an isomorphism of algebras provided
the product in
$P\otimes M$ is given by the formula
$$
(p\otimes x)(q\otimes y) = p\sw 1q\sw 1\otimes \Phi^{-1}(p\sw 2q\sw
2)\Phi(p\sw 3)x\Phi(q\sw 3)y,
$$
(cf.\ \cite[Theorem~11]{DoiTak:cle}).

In this case there is a left action of $P$ on $F$, $\mu: P\otimes
F\to F$  given
by $p\otimes f\mapsto \Phi(p f\su{-2})\Phi^{-1}(f\su{-1})f\sw 0 $.
When
$F$ is viewed as $P\otimes M$ this is simply the multiplication of
elements
in $P$, hence it is a left action. That fact that $\mu$ is compatible
with the
coactions is clear once one identifies $F$ with $P\otimes M$ via
$\Theta$. To check this
explicitly, take any   $f\in F$ and $p\in P$ and compute
\begin{eqnarray*}
\lambda(\mu(p\otimes f)) &=&
\lambda( \Phi(p f\su{-2})\Phi^{-1}(f\su{-1})f\sw 0) \\
&=&
p\sw 1f\su{-5}(Sf\su{-2})f\su{-1}\otimes
\Phi(p\sw 2 f\su{-4})\Phi^{-1}(f\su{-3})f\sw 0\\
& =& p\sw 1f\su{-2} \otimes \Phi(p f\su{-2})\Phi^{-1}(f\su{-1})f\su
0\\
&=& p\sw 1f\su{-1}\otimes \mu(p\sw 2\otimes f\su 0),
\end{eqnarray*}
as required. Thus if $F$ is a cleft extension of $M$, then
$H_{c}^{n}(P;F)=0$ for $n\ge 1$
and $H_{c}^{0}(P;F)\cong {}^{coP} F=M$.
\end{example}

For another example of a left $P$-comodule which is a left $P$-Hopf
module,
so that  the assumptions of Proposition~\ref{yooy} hold, take $F$ to
be
equal a $P$-bimodule $\Omega^nP$ of differential $n$-forms on $P$ in
a
left covariant differential structure (cf.\ Corollary~\ref{jjuu}
below).

\medskip
It will be convenient later to have an explicit description of this
cohomology
theory which does not involve $P$-invariance, so we give an
alternative formulation.

\begin{propos}\label{altfed}
       For a left $P$-comodule $F$ with coaction
$\lambda(f)=f_{[-1]}\tens f_{[0]}$, define a cochain complex
$(G^{*},\bar \extd)$ by $G^{n}=P^{n}\tens F$ for $n\ge 0$ with
derivation $\bar \extd f=1_{P}\tens f-\lambda(f)$ and
\begin{eqnarray*}
       \bar\extd(p_{1}\tens\dots\tens p_{n}\tens f) &=& 1_{P}\tens
p_{1}\tens\dots\tens p_{n}\tens f\,-\, \Delta(p_{1})\tens\dots\tens
p_{n}\tens f\,+\,\dots \cr &&+\,(-1)^{n}\,p_{1}\tens\dots\tens
\Delta(p_{n})\tens f\,-\, (-1)^{n}\,p_{1}\tens\dots\tens p_{n}\tens
\lambda(f)\ .
\end{eqnarray*}
Then there is a cochain isomorphism $\theta:(G^{*},\bar\extd)\to
       ({}^{coP}D^{*},\extd)$ given by $\theta(f)=S(f_{[-1]})\tens
f_{[0]}$
       and
       \begin{eqnarray*}
       \theta(p_{1}\tens\dots\tens p_{n}\tens f) &=&
       S(p_{1(1)})\tens p_{1(2)}\, S(p_{2(1)})\tens\dots\tens
       p_{n-1(2)}\, S(p_{n(1)}) \cr && \tens p_{n(2)}\,
S(f_{[-1]})\tens
f_{[0]}\ .
       \end{eqnarray*}
\end{propos}
\proof By explicit calculation the image of $\theta$ is in the
$P$-invariant part of $D^{n}$.  To show that $\theta$ is a cochain map
we calculate
\begin{eqnarray*}
       \theta(1_{P}\tens
p_{1}\tens\dots\tens p_{n}\tens f) &=& 1_{P}\tens S(p_{1(1)})\tens
p_{1(2)}\,
S(p_{2(1)})\tens\dots  \tens p_{n(2)}\, S(f_{[-1]})\tens f_{[0]}\ ,\cr
\theta(\Delta(p_{1})\tens\dots\tens
p_{n}\tens f) &=& S(p_{1(1)})\tens 1_{P}\tens p_{1(2)}\,
S(p_{2(1)})\tens\dots \tens p_{n(2)}\, S(f_{[-1]})\tens f_{[0]}\ ,\cr
\theta(p_{1}\tens\dots\tens
\Delta(p_{n})\tens f) &=& S(p_{1(1)})\tens p_{1(2)}\,
S(p_{2(1)})\tens\dots \tens 1_{P} \tens p_{n(2)}\, S(f_{[-1]})\tens
f_{[0]}\ ,\cr
\theta(p_{1}\tens\dots\tens p_{n}\tens
\lambda(f)) &=& S(p_{1(1)})\tens p_{1(2)}\,
S(p_{2(1)})\tens\dots \tens p_{n(2)}\, S(f_{[-1]})\tens 1_{P} \tens
f_{[0]}\ .
\end{eqnarray*}
Adding these terms with the appropriate signs shows that
$\theta\circ\bar\extd =
\extd\circ\theta$.  The inverse function is given by
\begin{eqnarray*}
       \theta^{-1}(p_{0}\tens f) &=& \eps(p_{0})\, f\ ,\cr
       \theta^{-1}(p_{0}\tens p_{1}\tens f) &=& \eps(p_{0})\, p_{1}\,
f_{[-1]}
       \tens f_{[0]}\ ,\cr
       \theta^{-1}(p_{0}\tens p_{1}\tens p_{2}\tens f) &=&
\eps(p_{0})\,
       p_{1}\,p_{2(1)}\, f_{[-2]} \tens p_{2(2)}\, f_{[-1]} \tens
f_{[0]}\ ,\cr
       \theta^{-1}(p_{0}\tens p_{1}\tens p_{2}\tens p_{3}\tens f) &=&
\eps(p_{0})\,
       p_{1}\,p_{2(1)}\,p_{3(1)}\, f_{[-3]} \tens p_{2(2)}\,p_{3(2)}
       \, f_{[-2]} \tens p_{3(3)} \, f_{[-1]} \tens f_{[0]}\ ,
\end{eqnarray*}
and so on.\quad$\square$

\section{Spectral sequences}
This section contains well known material on spectral sequences and
double complexes,
which we have taken from \cite{spseq}, though we have slightly
specialised the results.

\begin{defin}\rm
    A double complex of bidegree $(r,t)$ is a collection of vector
spaces $E^{n,m}$
(taken to be zero if either
$n<0$ or $m<0$) and linear maps $D:E^{n,m}\to E^{n+r,m+t}$ with
$D\circ
D=0$. A spectral sequence of degree $\ge s$ is a collection
$(E^{*,*}_r,D_r)$
of double complexes of bidegree $(r,1-r)$ for $r\ge s$ so that
\[
E^{n,m}_{r+1}\,\cong\,H^{n,m}(E_r,D_r)\,=\,(\ker D_r:E^{n,m}_r\to
E^{n+r,m+1-r}_r)/({\rm im}\,
D_r:E^{n-r,m+r-1}_r\to E^{n,m}_r)\ .
\]
\end{defin}

The problem is that knowledge of $D_r$ does not necessarily imply
knowledge of $D_{r+1}$,
though additional information, such as a ring structure, may help in
specific cases.
    The name of the game is to be
able to calculate the limit $\lim_{r\to\infty}E^{n,m}_{r}$. The
conditions that we put on the
vanishing of $E^{n,m}_{r}$ ensure that the $E^{n,m}_{r}$ stabilise
for sufficiently large $r$.

Next we give an example where spectral sequences arise, the example
that is essential
for the van Est spectral sequence.

\begin{example}\label{dcomp}\rm
Start with a double complex $C^{n,m}$ (taken to be zero if either
$n<0$ or $m<0$) with two differentials,
$\extd'$ of bidegree $(1,0)$ and $\extd''$ of bidegree $(0,1)$, which
satisfy
    $\extd''\circ \extd'+\extd'\circ \extd''=0$. The total complex
$T^s$
is defined as $T^s=\oplus_i C^{i,s-i}$ with differential
$\extd=\extd'+\extd''$.
Define
\begin{eqnarray*}
H_I^{n,m}(C) &=& H^{n,m}(C,d')\,=\,(\ker \extd':C^{n,m}\to
C^{n+1,m})/({\rm im}\,
\extd':C^{n-1,m}\to C^{n,m})\ ,\cr
H_{II}^{n,m}(C) &=& H^{n,m}(C,d'')\,=\,(\ker \extd'':C^{n,m}\to
C^{n,m+1})/({\rm im}\,
\extd'':C^{n,m-1}\to C^{n,m})\ .
\end{eqnarray*}
The map $\extd''$ induces a differential $\bar\extd''$ on the
bigraded complex $H_I^{n,m}(C)$
by $\bar\extd''[c]=[\extd'' c]$, and $\extd'$
induces a differential $\bar\extd'$ on the double complex
$H_{II}^{n,m}(C)$ by $\bar\extd'[c]=[\extd' c]$, where $[c]$
denotes the equivalence class of $c\in C^{n,m}$ under the quotient.
Now there are two spectral sequences $({}_I E^{*,*}_r,{}_I D_r)$ and
$({}_{II} E^{*,*}_r,{}_{II} D_r)$
for $r\ge 2$ with ${}_I E^{n,m}_2\cong
H^{n,m}(H_{II}^{*,*}(C),\bar\extd')$ and
${}_{II} E^{n,m}_2\cong H^{n,m}(H_{I}^{*,*}(C),\bar\extd'')$ which
both converge to the
the cohomology of the total complex $H^*(T,\extd)$.
This means that
$\oplus_{i} E^{i,n-i}_\infty=H^n(T,\extd)$, denoting the limits by
$E^{n,m}_\infty$.
\end{example}

\section{The van Est spectral sequence}
Again, $P$ is a Hopf algebra. The construction of the van Est
spectral sequence for Hopf algebras is based on the  following
example of a double complex.

\begin{example}\label{veexam}
\rm Suppose that there are left $P$-comodules $F^n$ for $n\ge 0$, and
left $P$-comodule maps $\bar\extd:F^n\to F^{n+1}$ with the property
that
$\bar \extd\circ\bar\extd=0$. Construct a double complex
$C^{n,m}=P^{\tens n}\tens F^m$ ($n,m\ge 0$) with differentials
$\extd':P^{\tens n}\tens F^m\to P^{\tens n+1}\tens F^m$ being the
Hopf-cochain
differentials given in Proposition~\ref{altfed} and
$\extd'':P^{\tens n}\tens F^m\to P^{\tens n}\tens F^{m+1}$ being
$(-1)^n\, \id\tens\bar\extd$. Note that the $(-1)^n$ factor is
included to
force the condition $\extd''\circ \extd'+\extd'\circ \extd''=0$.
\end{example}

\begin{lemma}
Suppose that the cochain complex $(F^n,\bar \extd)$ in
Example~\ref{veexam}
satisfies the additional
condition that each $F^n$ has a left $P$-action $:P\tens F^n\to F^n$
that
is a left $P$-comodule map, i.e., each of the $F^n$ is a left
$P$-Hopf module.
Then the spectral sequence ${}_I
E^{n,m}_r$ described
in Example~\ref{dcomp} converges to $H^s({}^{coP} F,\bar\extd)$.
\end{lemma}
\proof From Proposition~\ref{yooy},
$H_I^{n,m}(C)=H^{n,m}(C,\extd')=0$ for $n>0$
and $H_I^{0,m}(C)\cong {}^{coP} F^m$,
with the isomorphism mapping $f\in {}^{coP} F^m$ to $[1_P\tens f]\in
H_I^{0,m}(C)$. By definition
$\bar\extd''[1_P\tens f]=[1_P\tens \bar\extd f]$, so the isomorphism
identifies the complexes $(H_I^{0,m}(C),\bar\extd'')$ and $({}^{coP}
F^m,\bar \extd)$. Thus we have
${}_{II} E_2^{n,m}=  0$ for $n>0$ and   ${}_{II}
E_2^{0,m}=H^m({}^{coP}
F,\bar \extd) $. But now every
map ${}_{II}D_r$ ($r\ge 2$) must be zero, as it either maps into or
out of zero. This means that the
${}_{II}E_r^{*,*}$ spectral sequence stabilises at $r=2$, giving the
result.\quad$\square$

\begin{theorem}\label{ppuu} Suppose that $P$ is a Hopf algebra,
and that there is a cochain complex $(F^n,\bar \extd)$ of left
$P$-comodules
with the differential $\bar\extd$ being a comodule map.  Additionally
suppose that each $F^n$ has a left $P$-action $:P\tens F^n\to F^n$
which is a left $P$-comodule map.  Then there is a spectral sequence
beginning with $E^{n,m}_2=H^{n}_{c}(P;H^{m}(F^{*},\bar\extd))$ which
converges to $H^s({}^{coP} F,\bar\extd)$.
\end{theorem}
\proof We identify the spectral sequence ${}_{I} E^{n,m}_2$ in
Example~\ref{veexam}.
    First $\ker (\id_{P^{\tens n}}\tens\bar d)=P^{\tens n}\tens\ker
\bar
    \extd$ and ${\rm im}\, (\id_{P^{\tens n}}\tens\bar d)=P^{\tens
    n}\tens{\rm im}\, \bar \extd$.  The short exact sequence
\[
0 \longrightarrow ({\rm im}\, \bar\extd:F^{m-1}\to F^{m})
    \longrightarrow (\ker \bar\extd:F^m\to F^{m+1})
\longrightarrow H^m(F^*,\bar\extd) \longrightarrow 0
\]
remains exact if we tensor on the left with $P^{\tens n}$. Putting
these results together,
\[
H^m(P^{\tens n}\tens F^*,\id\tens\bar\extd)\,=\, P^{\tens n}\tens
H^m(F^*,\bar\extd)\ .
\]
Applying the induced differential to this gives the result.
\quad$\square$

\begin{cor} \label{jjuu}
       Suppose that $P$ is a Hopf algebra with bicovariant differential
calculus.
Then there is a spectral sequence beginning with
$E^{n,m}_2=H^{n}_{c}(P;H_{dR}^{m}(P))$ which converges to
$H^s({}^{coP}(\Omega^* P),\extd)$.
\end{cor}
\proof Put $(F^*,\bar\extd)$ equal to the de Rham complex $(\Omega^*
P,\extd)$ in Theorem~\ref{ppuu}.
The action of left multiplication $:P\tens \Omega^n P\to \Omega^n P$
is a
left $P$-comodule map. \quad$\square$

\begin{example}\rm
Suppose that $F$ is a cleft extension of $M$ as described in
Example~\ref{ex.cleft}. We thus
know that $F$ is a left $P$-Hopf module. If, in addition, $\Phi:P\to
F$ is an algebra map, then
the left action $\mu$ of $P$ on $F$ is induced from the product in
$F$ via the map $\Phi$, i.e.,
$\mu(p\ot f) = \Phi(p)f$. In this case, for any left $P$-covariant
differential structure $\Omega F$,
the $F$-bimodule of $n$-forms $\Omega^nF$ is a left $P$-module via
the map $\Phi$, i.e., there
are actions $\mu_n: P\otimes \Omega^nF \to \Omega^nF$, given by
$p\otimes \omega =
\Phi(p)\cdot \omega$. Since $\Omega F$ is a $P$-covariant calculus,
we can compute, for
any $p\in P$ and $\omega\in  \Omega^nF$,
$$
\lambda^n(\mu_n(p\ot \omega)) = \lambda^n(\Phi(p)\omega) = p\sw 1
\omega\su{-1}\ot \Phi(p\sw 2)\omega\su 0 = p\sw 1\omega\su{-1}\ot
\mu_n(p\sw 2\ot \omega),
$$
where $\lambda^n: \Omega^nF\to P\ot \Omega^nF$ is the $P$-coaction.
Thus each of the
$\Omega^nF$ is a left Hopf $P$-module, and we can infer from
Theorem~\ref{ppuu} that
there is a spectral sequence
beginning with $E^{n,m}_2=H^{n}_{c}(P;H^{m}(\Omega F,\bar\extd))$ that
converges to $H^s(M,\bar\extd)$. General constructions and explicit
examples of
differential structure on cleft extensions can be found in
\cite{BrzMaj:dif} (beware
  that \cite{BrzMaj:dif} uses right coactions rather than left
coactions).
\end{example}

\section{Braidings and $n$-forms}\label{uiui}
The material contained in this section has existed for a long time,
and was mostly developed by Woronowicz in \cite{worondiff}. Thus we
merely give a brief
description and point out that
proofs can be seen in \cite{worondiff}.

For a Hopf algebra $P$ with a bicovariant differential structure,
there are right
and left coactions $\rho$ and $\lambda$, given by
\[
\rho(p.\extd q)\,=\,p_{(1)}.\extd q_{(1)} \tens p_{(2)}\, q_{(2)}\
,\quad
\lambda(p.\extd q)\,=\,p_{(1)}\, q_{(1)} \tens p_{(2)}.\extd q_{(2)}\
,
\]
for all $p,q\in P$. We shall use the index notation
$\rho(\xi)=\xi_{[0]}\tens\xi_{[1]}$ and
$\lambda(\xi)=\xi_{[-1]}\tens\xi_{[0]}$ (summation understood).
$\Omega^1P$ is a bicomodule with these coactions, meaning that the
left and right
coactions commute, i.e.,
we can write the following without ambiguity:
\[
(\lambda\tens\id)\rho(\xi)\,=\,(\id\tens\rho)\lambda(\xi)\,=\,
\xi_{[-1]}\tens\xi_{[0]}\tens\xi_{[1]}\ .
\]
Denote by $L^1$ the left invariant 1-forms ${}^{coP}(\Omega^1 P)$.
There is a map $Y:\Omega^1 P\to P\tens L^1$ defined by
$Y(\xi)=\xi_{[-2]}\tens S(\xi_{[-1]}).
\xi_{[0]}$. In terms of $P$-actions and coactions we get
$Y:{}^\bullet_\bullet\Omega^1 P^\bullet_\bullet \to {}^\bullet_\bullet
P^\bullet_\bullet \tens (L^1)^\bullet_\bullet $. The upper (resp.\
lower) dots indicate what coaction (resp.\ action) we take, i.e., the
right action
    and coaction are the tensor product ones, wheras
    the left action and coaction are purely on the first component.
The right action on $L^1$ is $\chi\ra p=S(p{(1)}).\chi.p{(2)}$.
Note that as a right $P$-module and a right $P$-comodule $L^1$
satisfies the
Yetter-Drinfeld condition. The main consequences of this are
explained in the following remarks.

\begin{remark}\label{yydd} \rm A tensor category consisting of
objects which are
right $P$-modules and right $P$-comodules, with morphisms which are
right $P$-module and right $P$-comodule maps, is said to satisfy the
(right)
{\em Yetter-Drinfeld condition} if $\rho(\eta\ra a)=\eta_{[0]} \ra
a_{(2)}
\tens S(a_{(1)})\,\eta_{[1]}\,a_{(3)}$.  Here $\eta$ is any element
of any object $V$ and $a$ is an element of $P$. The symbol $\ra$
indicates the right action of $P$ on $V$ and $\eta_{[0]} \tens
\eta_{[1]}$ denotes the right coaction. The tensor product has the
usual tensor product action and coaction, and the associator is
trivial. Objects of a category that satisfies the (right)
Yetter-Drinfeld condition are known as {\em Yetter-Drinfeld} or {\em
crossed} modules. The map $\sigma_{VW}:V\tens W\to W\tens V$ defined
by $
\sigma(\xi\tens \eta)=\eta_{[0]}\tens \xi\ra \eta_{[1]} $ is a
braiding for the category.  If the antipode $S$ is invertible, there
is an inverse braiding $ \sigma^{-1}(\eta\tens\xi)=\xi\ra
S^{-1}(\eta_{[1]})\tens \eta_{[0]} $.

We call an element $x\in V\tens V$ {\em symmetric} provided
$\sigma(x)=x$, and define
$V\wedge V$ to be the quotient of $V\tens V$ by the subspace $\cS(V)$
of symmetric
elements. In the same manner we can define
$V\wedge V\wedge V$ as the quotient of $V\tens V\tens V$ by the
subspace
generated by $\cS(V)\tens V$ and $V\tens \cS(V)$, and so on.
\end{remark}

\begin{remark} \label{remser} \rm We shall take the isomorphism $Y$
mentioned earlier seriously and
{\em define} ${}^\bullet_\bullet\Omega^1 P^\bullet_\bullet =
{}^\bullet_\bullet
P^\bullet_\bullet \tens (L^1)^\bullet_\bullet$. This naturally leads
to the definitions
${}^\bullet_\bullet\Omega^2 P^\bullet_\bullet={}^\bullet_\bullet
P^\bullet_\bullet\tens (L^1)^\bullet_\bullet\wedge
(L^1)^\bullet_\bullet$,
${}^\bullet_\bullet\Omega^3 P^\bullet_\bullet={}^\bullet_\bullet
P^\bullet_\bullet\tens (L^1)^\bullet_\bullet\wedge
(L^1)^\bullet_\bullet\wedge
(L^1)^\bullet_\bullet$ etc. The wedge product making $\Omega^* P$
into a graded algebra is
    $(p\tens v)
\wedge (q\tens w)=p\,q_{(1)}\tens v\ra q_{(2)}\wedge w$, for all
$p,q\in P$ and $v,w\in L^1$. The de Rham
differential on the
left invariant 1-forms is given by
\begin{eqnarray*}
\extd(S(p_{(1)}).\extd p_{(2)}) &=& \extd S(p_{(1)})\wedge\extd
p_{(2)} \,=\, -\,S(p_{(1)}).\extd p_{(2)}
\wedge S(p_{(3)}).\extd p_{(4)} \ ,
\end{eqnarray*}
and by explicit calculation this $\extd:L^1 \to L^1\wedge L^1$ is a
right comodule map. To derive this equation we used the result
$$0 = \extd(\eps(p)) = \extd(S(p\sw 1)p\sw 2) = \extd(S(p\sw 1)).p\sw
2 - S(p\sw 1).\extd(p\sw 2),
$$
so that
$\extd(S(p)) = -S(p\sw 1).\extd p\sw 2. S(p\sw 3)$.
The differential is extended to $\extd:(L^1)^{\wedge n} \to
(L^1)^{\wedge n+1}$ by
\begin{eqnarray*}
\extd(\xi_1\wedge\xi_2\wedge\dots\wedge\xi_n)
&=&\extd\xi_1\wedge\xi_2\wedge\dots\wedge\xi_n
\,-\, \xi_1\wedge\extd\xi_2\wedge\dots\wedge\xi_n \cr
&& +\, (-1)^{n+1}\, \xi_1\wedge\xi_2\wedge\dots\wedge\extd\xi_n\ ,
\end{eqnarray*}
and to $\extd:\Omega^n P\to\Omega^{n+1}P$ by
$\extd(p\tens v)=p_{(1)}\tens S(p_{(2)}).\extd p_{(3)} \wedge
v+p\tens\extd v$.
\end{remark}

\begin{remark}\rm
The discussion of the preceding remarks allows us in principle to
calculate the cohomology of the left
invariant forms
$(L^1)^{\wedge n}$, given the (often finite dimensional)
right $P$ module and comodule $L^1$ and $\extd:L^1 \to L^1\wedge L^1$.
However this construction does not use anything corresponding to the
Lie algebra of
a Lie group.

Incidently, though the braiding introduced here may seem rather
arbitrary, it is not too
difficult to justify. One way is to see that it is the braiding (in
the sense of \cite{Madore})
corresponding
    to the left covariant derivative on the bimodule $\Omega^1 P$
which kills all left invariant forms. There will be a better reason
later.

\end{remark}

\section{Adjoint coactions and the  Hopf-Lie algebra}
In the case of a Lie group $G$, the left Adjoint action $:G\times
G\to G$
is given by $\Ad_g(h)=ghg^{-1}$. Differentiating this in the second
variable
$h$ in the direction of $v$ in the Lie algebra $\gg$, we get
$\Ad_{g*}(v)=gvg^{-1}\in\gg$.  Finally differentiating with respect to
$g$ in the direction $w\in\gg$ we get the Lie bracket $[w,v]=wv-vw$.
By following this prescription for a Hopf algebra we get the most
direct justification for the braiding in section \ref{uiui}.

Again $P$ is a Hopf algebra with a bicovariant differential calculus
$\Omega^*P$.

\begin{defin}\rm
For a Hopf algebra $P$, the {\em Hopf-Lie algebra} is defined as
\[
\gp\,=\,\{\psi:\Omega^1 P\to
k:\psi(\xi.p)=\psi(\xi)\,\eps(p)\quad\forall p\in P\}\ .
\]
The left adjoint $P$-coaction on $P$ is defined by
$\Ad^L(p)=p_{(1)}\,S(p_{(3)})\tens p_{(2)}$. With the help of this
coaction we define a {\em bracket},
    for all $\alpha,\beta\in
\gp$,
\[
[\alpha,\beta]\,=\,\alpha\circ \extd\circ(\id\tens\beta)\circ
\Pi_0\circ \Ad^L_*:\Omega^1 P\to k\ .
\]
\end{defin}

This definition of a Hopf-Lie algebra has a classical motivation.
Classically, the Lie algebra of a Lie group can be identified with a
space dual to the cotangent space. In the case of a general Hopf
algebra $P$, $k$ is a right $P$-module with the action given by the
counit $\eps$ (remember that $\eps$ is a character of $P$). Also
$\Omega^1P$ is a right $P$-module, and so $\gp$ is simply a space of
right $P$-linear maps $\Omega^1P\to k$.

\begin{propos} \label{oouu} \cite{worondiff} The Hopf-Lie algebra is
closed under
the bracket $[\bullet , \bullet ]$, i.e., for all $\alpha,\beta\in
\gp$,  $[\alpha,\beta]:\Omega^1 P\to k$
is in $\gp$. Furthermore
\[
[\alpha,\beta](\xi) \,=\, \alpha(\extd(\xi_{[-1]}))\, \beta(\xi_{[0]})
- \alpha(\xi_{[-1]}\,S(\xi_{[1]}).\extd(\xi_{[2]}))\,
\beta(\xi_{[0]})\ .
\]
\end{propos}
\proof To show that $[\alpha,\beta]\in\gp$, using the fact that
$\beta\in\gp$, we write
\begin{eqnarray}
(\id\tens\beta)\circ \Pi_2\circ \Ad^L_*(\xi.p) &=& (\id\tens\beta)
(\xi_{[-1]}\,p_{(1)}\,
S(p_{(3)})\,S(\xi_{[1]})\tens
\xi_{[0]}\,p_{(2)}) \cr
&=& \xi_{[-1]}\,p_{(1)}\,
S(p_{(3)})\,S(\xi_{[1]})\,
\beta(\xi_{[0]}\,p_{(2)})\cr
&=& \xi_{[-1]}\,p_{(1)}\,
S(p_{(3)})\,S(\xi_{[1]})\,
\beta(\xi_{[0]})\,\eps(p_{(2)})\cr
&=& \xi_{[-1]}\,p_{(1)}\,
S(p_{(2)})\,S(\xi_{[1]})\,
\beta(\xi_{[0]})\cr
&=& \xi_{[-1]}\,S(\xi_{[1]})\,
\beta(\xi_{[0]})\,\eps(p)\ .
\end{eqnarray}
Next we have, using the fact that $\extd$ is a derivation,
\begin{eqnarray}
[\alpha,\beta](\xi) &=& \alpha(\extd(\xi_{[-1]}\,S(\xi_{[1]})))\,
\beta(\xi_{[0]}) \cr
&=& \alpha(\extd(\xi_{[-1]}).S(\xi_{[1]}))\,
\beta(\xi_{[0]}) + \alpha(\xi_{[-1]}.\extd(S(\xi_{[1]})))\,
\beta(\xi_{[0]})
\cr
&=& \alpha(\extd(\xi_{[-1]}))\,\eps(S(\xi_{[1]}))\,
\beta(\xi_{[0]}) +
\alpha(\xi_{[-1]}.\extd(S(\xi_{[1]})))\,\eps(\xi_{[2]})\,
\beta(\xi_{[0]})\cr
&=& \alpha(\extd(\xi_{[-1]}))\,
\beta(\xi_{[0]}) + \alpha(\xi_{[-1]}.\extd(S(\xi_{[1]})).\xi_{[2]})\,
\beta(\xi_{[0]})\cr
&=& \alpha(\extd(\xi_{[-1]}))\,
\beta(\xi_{[0]}) - \alpha(\xi_{[-1]}\,S(\xi_{[1]}).\extd(\xi_{[2]}))\,
\beta(\xi_{[0]})\ .\quad\square
\end{eqnarray}

The next result justifies the definition of $\extd$ on the
left invariant 1-forms given in Remark~\ref{remser}, and shows that
$\extd$ on the
left invariant 1-forms is dual to the Lie bracket on $\gp$. The
$\id\tens\id-\sigma$ appearing in the formula allows us to take the
quotient from the tensor product to the wedge product.

\begin{propos} \label{coyy} \cite{worondiff} For any $p\in P$, let
$\xi=S(p_{(1)}).\extd p_{(2)}$ be the corresponding  left invariant
1-form. Then
\begin{eqnarray*}
[\alpha,\beta](\xi)
    &=&
\ev(\id\tens\ev\tens\id)((\alpha\tens\beta)\tens(\id\tens\id-\sigma)(-\,
S(p_{(1)}).\extd p_{(2)} \tens S(p_{(3)}).\extd p_{(4)}))\\
&=&\ev(\id\tens\ev\tens\id)((\alpha\tens\beta)\tens(\id\tens\id-\sigma)
(\extd\xi).
\end{eqnarray*}
\end{propos}
\proof
This is proven by explicit calculation, beginning from \ref{oouu}.
\begin{eqnarray*}
[\alpha,\beta](\xi) &=&  - \,\alpha(S(\xi_{[1]}).\extd(\xi_{[2]}))\,
\beta(\xi_{[0]}) \cr
&=& - \,\alpha(S(S(p_{(2)})\, p_{(5)}).\extd(S(p_{(1)})\, p_{(6)}))\,
\beta(S(p_{(3)}).\extd p_{(4)})
\cr
&=& - \,\alpha(S(p_{(5)})\,S^2(p_{(2)})\,S(p_{(1)}).\extd p_{(6)})\,
\beta(S(p_{(3)}).\extd p_{(4)}) \cr
&& - \,\alpha(S(p_{(5)})\,S^2(p_{(2)}).\extd S(p_{(1)}). p_{(6)})\,
\beta(S(p_{(3)}).\extd p_{(4)})\cr
&=& - \,\alpha(S(p_{(3)}).\extd p_{(4)})\,
\beta(S(p_{(1)}).\extd p_{(2)}) \cr
&& - \,\alpha(S(p_{(5)})\,S^2(p_{(2)}).\extd S(p_{(1)}))\,
\beta(S(p_{(3)}).\extd p_{(4)})\cr
&=& - \,\alpha(S(p_{(3)}).\extd p_{(4)})\,
\beta(S(p_{(1)}).\extd p_{(2)}) \cr
&& + \,\alpha(S(p_{(7)})\,S^2(p_{(4)})\,S(p_{(1)}).\extd p_{(2)}.
S(p_{(3)}))\,
\beta(S(p_{(5)}).\extd p_{(6)})\cr
&=& - \,\alpha(S(p_{(3)}).\extd p_{(4)})\,
\beta(S(p_{(1)}).\extd p_{(2)}) \cr
&& + \,\alpha(S(p_{(6)})\,S^2(p_{(3)})\,S(p_{(1)}).\extd p_{(2)})\,
\beta(S(p_{(4)}).\extd p_{(5)})\cr
&=& - \,\alpha(S(p_{(3)}).\extd p_{(4)})\,
\beta(S(p_{(1)}).\extd p_{(2)}) \cr
&& + \,\alpha((S(p_{(1)}).\extd p_{(2)}) \ra S(p_{(3)})\,p_{(6)})\,
\beta(S(p_{(4)}).\extd p_{(5)}) \cr
&=& -\ev(\id\tens\ev\tens\id)((\alpha\tens\beta)\tens(
S(p_{(1)}).\extd p_{(2)} \tens S(p_{(3)}).\extd p_{(4)} \cr
&&-\,
S(p_{(4)}).\extd p_{(5)} \tens (S(p_{(1)}).\extd p_{(2)}) \ra
S(p_{(3)})\,p_{(6)})),
\end{eqnarray*}
as required for the first equality.
   The second equality follows from the formula for $\extd$ in
\ref{remser}. \quad
$\square$\medskip

\section{The Hopf-Lie algebra cohomology}
Here we continue from the last section, and ask what the analogue
of the Lie algebra cohomology is for a Hopf algebra with bicovariant
differential calculus. Without the Jacobi identity the standard
formula will not work, we have to include the braiding
in the definition of the cochain complex. This section will
show how this can be done.

To define the Lie algebra cohomology, we need to transfer the
braiding on the
left invariant forms to the Lie algebra. There is an evaluation map
$:\gp\tens L^1\to k$, and this identifies the Lie algebra with the
dual of the
left invariant forms. Now we look at the dual operation on objects of
a Yetter Drinfeld category.

\begin{remark} \label{YDdual} \rm Suppose that the antipode $S$ on
$P$ is
invertible.  Following on from Remark~\ref{yydd}, we define the dual
of an
object $V$ in a  category with the Yetter-Drinfeld condition to be
the vector space dual
$V^*$ with right action and coaction
\[
(\alpha\ra p)(v)\,=\,\alpha(v\ra S^{-1}(p))\ ,\quad
\alpha_{[0]}(v)\, \alpha_{[1]}\,=\,\alpha(v_{[0]})\, S(v_{[1]})\ ,
\]
for all $p\in P$, $v\in V$ and $\alpha\in V^*$. The evaluation map
$\ev:V^*\tens V\to k$ preserves the action and coaction.
As the action and coaction on $V^{*}$ satisfy the Yetter-Drinfeld
condition
we can define a braiding $\sigma:V^{*}\tens V^{*}\to V^{*}\tens V^{*}$
by the usual formula.
The braiding on $V$ and $V^*$ are connected by the following formula,
for $\alpha,\beta\in V^*$ and $v,w\in V$:
\begin{eqnarray}\label{olol}
\ev(\id\tens\ev\tens\id)(\alpha\tens\beta\tens\sigma(v\tens w))\,=\,
\ev(\id\tens\ev\tens\id)(\sigma(\alpha\tens\beta)\tens v\tens w)\ .
\end{eqnarray}
\end{remark}

\begin{remark} \rm By using equation~(\ref{olol}) we can identify
$(L^{1})^{\wedge
n}$ with the linear maps from $\gp^{\wedge n}$ to $k$.  To avoid
unnecessary crossings in the braided category we use the evaluation
which pairs off the elements from the inside, i.e.\
$\ev(\id\tens\ev\tens\id):
\gp\wedge\gp\tens L^{1}\wedge L^{1}\to k$ etc., which we just write
$\underline{\ev}$.

If $\gp$ is finite dimensional, the map
$\id\tens\id-\sigma:\gp\wedge\gp\to
\gp\wedge\gp$ is invertible, and we write $T=[,]\circ
(\id\tens\id-\sigma)^{{-1}}:\gp\wedge\gp \to\gp$.  For example
if $\sigma$ has eigenvalues $\pm 1$ (the commutative case) then
$T(\alpha\wedge\beta)=\frac12[\alpha,\beta]$.  In the general case,
Proposition~\ref{coyy} gives
\begin{eqnarray}\label{olol44}
       \underline{\ev}((\alpha\wedge\beta)\tens
\extd\xi)\,=\,
T(\alpha\wedge\beta)(\xi)\ .
\end{eqnarray}
We can continue this formula, e.g.\
\begin{eqnarray}\label{olol55}
       \underline{\ev}((\alpha\wedge\beta\wedge\gamma)\tens
       \extd(\xi\wedge\eta)) &=&
\underline{\ev}((\alpha\wedge\beta\wedge\gamma)\tens
       (\extd\xi\wedge\eta-\xi\wedge\extd\eta)) \cr
&=& \underline{\ev}(\alpha\wedge
T(\beta\wedge\gamma)-T(\alpha\wedge\beta)\wedge\gamma)(\xi\wedge\eta)\
.
\end{eqnarray}
\end{remark}

Now we are in position to introduce the cohomology of a Hopf-Lie
algebra.
\begin{defin}\rm
       Define a cochain complex by $K^{n}$ being
       the linear maps from $\gp^{\wedge n}$ to $k$.  The differential
       $\extd:K^{n}\to K^{n+1}$ is given by
       \begin{eqnarray}
       \extd\phi(\alpha_{0}\wedge\dots\wedge\alpha_{n}) &=& \phi(
       \alpha_{0}\wedge\dots\wedge
T(\alpha_{n-1}\wedge\alpha_{n}))-\dots \cr
       && -\, (-1)^{n}
       \phi(T(\alpha_{0}\wedge\alpha_{1})\wedge\dots\wedge\alpha_{n})
\ .
       \end{eqnarray}
       The cohomology $H_{HL}(\gp)$ of this complex is known as  the
{\em Hopf-Lie
cohomology of $\gp$}.
\end{defin}

\begin{propos}
The Hopf-Lie
cohomology $H_{HL}(\gp)$ of $\gp$ is isomorphic to
the cohomology $H^s({}^{coP}(\Omega^* P),\extd)$
of the left invariant forms.
\end{propos}
\proof Combining the results of this section and the last section.

\section*{Acknowledgements}
T.\ Brzezi\'nski  would like to thank the Engineering and Physical
Sciences Research
Council for an Advanced Fellowship.

\end{document}